# Fault Analysis Using Gegenbauer Multiresolution Analysis

Luciana R. Soares and Hélio M. de Oliveira, *Member, IEEE*

*Abstract*—This paper exploits the multiresolution analysis in the fault analysis on transmission lines. Faults were simulated using the ATP (Alternative Transient Program), considering signals at 128/cycle. A nonorthogonal multiresolution analysis was provided by Gegenbauer scaling and wavelet filters. In the cases where the signal reconstruction is not required, orthogonality may be immaterial. Gegenbauer filter banks are thereby offered in this paper as a tool for analyzing fault signals on transmission lines. Results are compared to those ones derived from a 4-coefficient Daubechies filter. The main advantages in favor of Gegenbauer filters are their smaller computational effort and their constant group delay, as they are symmetric filters.

## I. INTRODUCTION

Discrete wavelet transforms associated with a multiresolution analysis (MRA) have already been stated as a powerful tool for fault analysis on transmission lines [1]-[8]. The choice of an appropriate mother wavelet is the earliest step when performing such an analysis. Orthogonal filter banks (Daubechies, Symmlet, Coiflet, etc.) are often preferred for analyzing power system transients [1]-[8]. Nevertheless, this paper shows that signals derived from a nonorthogonal multiresolution analysis can alternatively provide effective fault detection and fault location algorithms.

Recently, new wavelet families were introduced [9]-[11] with the basis of the idea of linking 2$^{nd}$-order ordinary differential equations with the transfer function of multiresolution analysis filters. These families are related to Legendre, Chebyshev and Gegenbauer polynomials, the latter being a broad class that generalizes the earliest two. Each one of such wavelets has compact support, and they are furthermore derived from constant group delay filters.

Although Gegenbauer MRA throws away orthogonality, this may be irrelevant when the analysis does not require the reconstruction of the analyzed signals. Indeed, Gegenbauer filter banks seem not be recommended for signal compression, since signal reconstruction may furnish poor approximations of the original signal.

The Gegenbauer MRA is first revisited and an example of application (transmission line protection) is illustrated. Results derived from Gegenbauer filter banks are compared to the ones derived from a Daubechies filter bank. Gegenbauer filters are offered as a tool for analyzing fault signals on transmission lines. It can be an alternative to asymmetric orthogonal filters, especially when a fast computing is essential or when group delay plays an important role.

## II. GEGENBAUER MULTIRESOLUTION ANALYSIS

Gegenbauer polynomials are solution of the differential equation, *n* integer:

$$(1-z^2)\frac{d^2y}{dz^2} - (2\alpha+1)z\frac{dy}{dz} + n(n+2\alpha)y = 0. \quad (1)$$

The *n*$^{th}$-order orthogonal Gegenbauer polynomial $C_n^{(\alpha)}(z)$ can be found, for $n>2$, $|z|\leq 1$ and $\alpha > -1/2$, by the recurrence relation [12]:

$$C_n^{(\alpha)} = \frac{1}{n}\left[2(\alpha+n-1)zC_{n-1}^{(\alpha)}(z) - (2\alpha+n-2)C_{n-2}^{(\alpha)}(z)\right] \quad (2)$$

with $C_1^{(\alpha)}(z) = 2\alpha z$ and $C_2^{(\alpha)}(z) = 2\alpha(\alpha+1)z^2 - \alpha$. Additionally, $C_n^{(\alpha)}(-z) = (-1)^n C_n^{(\alpha)}(z)$ holds.

It was proved [11] that the variable change $z=\cos(\omega/2)$, $n=\nu$, $\nu$ odd, under the constraint $\alpha$ strictly positive, generates a low-pass frequency selective FIR filter, in such a way that its impulse response converges to a scaling function of a MRA. Although Gegenbauer polynomial holds orthogonality for $\alpha > -1/2$, it has not a low-pass behavior within the interval $-1/2 < \alpha < 0$ [11]. Furthermore, 1$^{st}$-order Chebyshev polynomials (any $\nu$$^{th}$-order Gegenbauer polynomial with $\alpha=0$) do not generate scaling functions of a MRA [10].

The Gegenbauer scaling filter was defined by [11]:

$$H_\nu(\omega) = \frac{C_\nu^{(\alpha)}(\cos(\omega/2))}{C_\nu^{(\alpha)}(1)} \cdot e^{-j\frac{\omega\nu}{2}}, \quad (3)$$

where $C_\nu^{(\alpha)}(.)$ is a $\nu$$^{th}$-order Gegenbauer polynomial.

After some trigonometric handling, it can be shown that the

This work was partially supported by Brazilian Ministry of Education (CAPES) and by the Brazilian National Council for Scientific and Technological Development (CNPq) under research grant #306180.
L. R. Soares and H. M. de Oliveira are with the Signal Processing Group, Department of Electronics and Systems, Federal University of Pernambuco – UFPE, C.P. 7800, 50711-970, Recife, PE, Brazil (e-mail: lusoares@ufpe.br, hmo@ufpe.br).

Gegenbauer filter coefficients are given by [11]:

$$\frac{h_k^\nu}{\sqrt{2}} = \frac{1}{C_\nu^{(\alpha)}(1)} \frac{\Gamma(\alpha+k)\Gamma(\alpha+\nu-k)}{k!(\nu-k)!\Gamma^2(\alpha)}, \quad (4)$$

where $k = 0,1,\ldots,\nu$, and $\Gamma(.)$ is the standard gamma function.

The number of zeroes within the interval $0 < \omega \leq 2\pi$ depends on the degree ($\nu$) of the Gegenbauer polynomial. Moreover, all zeroes in the Z-plane are located on the unit circle. For a preset $\nu$, the parameter $\alpha$ also rules the main lobe width as well as the stop-band attenuation of the Gegenbauer scaling filters. Illustrative examples of such scaling filters are shown in the Fig. 1. It can be seen that these scaling filters are low-pass frequency selective FIR filters.

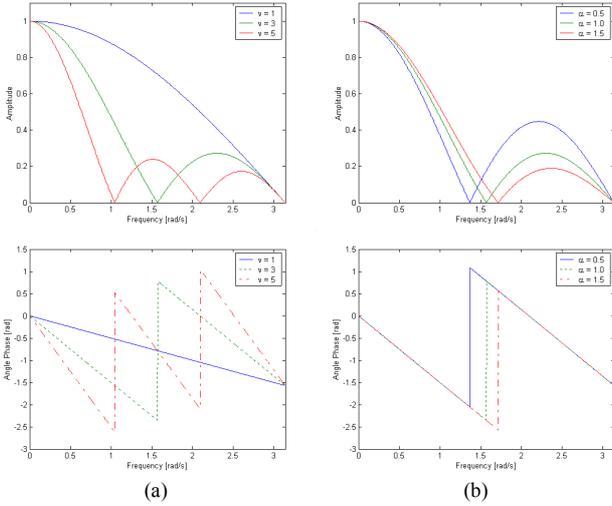

Fig. 1. Frequency response of Gegenbauer scaling filters: (a) $\alpha = 1$ and $\nu = 1, 3, 5$; (b) $\alpha = 0.5, 1.0, 1.5$ and $\nu = 3$.

Despite the fact of Gegenbauer scaling filters always hold the *necessary condition* for orthogonality, $\nu=1$ is the sole Gegenbauer filter order that generates an orthogonal MRA. Actually, the Gegenbauer filter ($\nu=1$, any $\alpha$) collapses into the Haar filter, which is the only orthogonal symmetrical filter [13].

A Gegenbauer wavelet filter of a nonorthogonal MRA can be implemented through a $\pi$-shift on the frequency response of the Gegenbauer scaling function. In discrete time domain, such a frequency-shift corresponds to a circular shift in the low-pass filter coefficients, i.e., $g_k = (-1)^k h_{1-k}$, $k=1,2,\ldots,\nu$, where $g_k$ are the wavelet filter coefficients.

Defining a nonorthogonal filter bank based on low-pass and high-pass Gegenbauer filters, the scaling and wavelet waveforms can be derived by an iterative procedure through its filter coefficients [14]. Fig. 2 shows a few Gegenbauer scaling and wavelet functions.

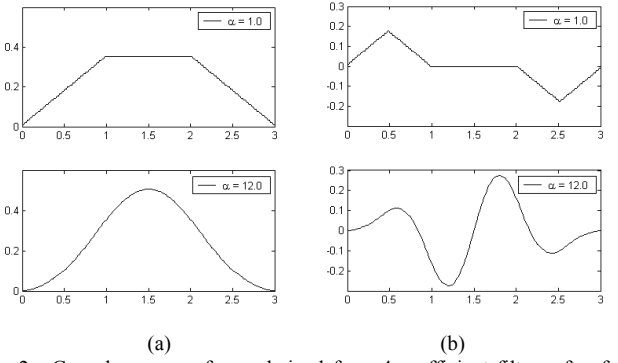

Fig. 2. Gegenbauer waveforms derived from 4-coefficient filters after four iterations. Parameters: $\nu = 3$ and $\alpha = 1.0, 12.0$. (a) scaling functions; (b) wavelet functions.

Since $\alpha$ can assume countless values ($\alpha \in \mathbf{R}$; $\alpha > 0$), the Gegenbauer MRA family has uncountable Gegenbauer scaling functions and mother wavelets for any $\nu$ (positive) odd. Let Geg$\nu a\alpha$ denote a MRA with parameters $\nu$ and $\alpha$, e.g., Geg3a1 denotes a Gegenbauer analysis with $\nu=3$ and $\alpha=1$. Special cases of the MRA Gegenbauer family are Haar ($\nu=1$, any $\alpha$), Legendre (any $\nu$, $\alpha=0.5$) [9], and 2nd-order Chebyshev family (any $\nu$, $\alpha=1.0$) [10].

Gegenbauer MRA family presents scaling filters with even symmetry and wavelet filters with odd symmetry. They are therefore, type II (even length and even symmetry) and type IV (even length and odd symmetry) FIR generalized linear phase filters [15], respectively. Every single one filter has linear phase and constant group delay, given by $\nu/2$, which means that *no different delay is introduced at different frequencies* of the analyzed signal. Another advantage in favor of symmetrical filters is concerning the computational effort, since only a half (or even lesser) of its filter coefficients needs to be computed.

### III. FAULT ANALYSIS

Improvements in fault analysis on transmission lines have been achieved via wavelet-based algorithms both for schemes that use the fundamental frequency components of the fault signals [4]-[7], and for those ones that handle with high-frequency components, being thereby based on the wave-traveling theory [8].

Considering merely single-ended algorithms for fault location and signals sampled at 128/cycle, methods based on the fundamental components might be more appropriated [16]. The wavelet-based fault-location algorithm presented in [5] provides a good estimate for the fault distance from the monitoring terminal.

*A. On the Fault Analysis Algorithm*

An estimation of the fault location can be derived by means of the apparent impedance approach [17]. Instead of dealing with the original voltages and currents, approximated versions of these signals are used to find out the fundamental components via a 1-cycle windowed Fourier transform.

Suppose that the line monitoring signals are sampled at 128/cycle during 8 cycles, and that the length $N$ of the signals is expressed by $N = 2^j$. Signals can thus be decomposed until a maximum level $j_{max}=10$. Table I shows that after the 3rd-level decomposition, approximated versions of the original signal are represented by less than 16 samples/cycle, which is not suitable for fault analysis [18]. Consequently, the 3rd-level is selected to derive the approximated versions of the original signals, which are to be used in the fault location algorithm. Table I also shows that the orthogonal 3rd-level MRA approximated versions have only frequencies below 480 Hz. The cut-off frequency was obtained at –3 dB.

TABLE I
FREQUENCY BANDS OF AN ORTHOGONAL MULTIRESOLUTION ANALYSIS AT DIFFERENT DECOMPOSITION LEVELS[a]

| Level ($j$) | Signal length (samples/cycle) | Scaling Filter (Hz) | Wavelet Filter (Hz) |
|---|---|---|---|
| 1 | 64 | 0 – 1920 | 1920 – 3840 |
| 2 | 32 | 0 – 960 | 960 – 1920 |
| 3 | 16 | 0 – 480 | 480 – 960 |
| 4 | 8 | 0 – 240 | 240 – 480 |
| 5 | 4 | 0 – 120 | 120 – 240 |
| 6 | 2 | 0 – 60 | 60 – 120 |
| 7 | 1 | 0 – 30 | 30 – 60 |

These frequency bands can be more selective when using the Gegenbauer MRA. Tables II and III show the frequency bands for 3rd and 7th-order filter, respectively. It can be seen from the table II that approximated versions at the 3rd-level have frequency components below 219 Hz and 278 Hz, considering $\alpha=1$ and $\alpha=12$, respectively. At table III, these frequencies are below 107 Hz and 171 Hz. The frequency values presented at tables II and III are rounded.

TABLE II
FREQUENCY BANDS AT SEVERAL DECOMPOSITION LEVELS FOR A FEW 3RD-ORDER GEGENBAUER MULTIRESOLUTION ANALYSIS[a]

| Level ($j$) | Scaling Filter (Hz) | | Wavelet Filter (Hz) | |
|---|---|---|---|---|
| | Geg3a1 | Geg3a12 | Geg3a1 | Geg3a12 |
| 1 | 0 – 877 | 0 – 1110 | 2963 – 3840 | 2730 – 3840 |
| 2 | 0 – 439 | 0 – 555 | 1481 – 1920 | 1365 – 1920 |
| 3 | 0 – 219 | 0 – 278 | 741 – 960 | 682 – 960 |
| 4 | 0 – 110 | 0 – 139 | 370 – 480 | 341 – 480 |
| 5 | 0 – 55 | 0 – 69 | 185 – 240 | 171 – 240 |
| 6 | 0 – 27 | 0 – 35 | 93 – 120 | 85 – 120 |
| 7 | 0 – 14 | 0 – 17 | 46 – 60 | 43 – 60 |

TABLE III
FREQUENCY BANDS AT SEVERAL DECOMPOSITION LEVELS FOR A FEW 7TH-ORDER GEGENBAUER MULTIRESOLUTION ANALYSIS[a]

| Level ($j$) | Scaling Filter (Hz) | | Wavelet Filter (Hz) | |
|---|---|---|---|---|
| | Geg7a1 | Geg7a12 | Geg7a1 | Geg7a12 |
| 1 | 0 – 427 | 0 – 683 | 3413 – 3840 | 3157 – 3840 |
| 2 | 0 – 214 | 0 – 341 | 1706 – 1920 | 1579 – 1920 |
| 3 | 0 – 107 | 0 – 171 | 853 – 960 | 789 – 960 |
| 4 | 0 – 53 | 0 – 85 | 427 – 480 | 395 – 480 |
| 5 | 0 – 27 | 0 – 43 | 213 – 240 | 197 – 240 |
| 6 | 0 – 13 | 0 – 21 | 107 – 120 | 99 – 120 |
| 7 | 0 – 7 | 0 – 11 | 53 – 60 | 49 – 60 |

Due to its frequency-selective response, scaling

---
[a] Signals sampled at 128 by cycle.

Gegenbauer filters may even be more fitting to fault location than scaling Daubechies filters. It is worthy to remark that frequency bands at the main lobe can be selected by a proper choice of $v$ and $\alpha$.

Because of the relationship between the fundamental components of voltage and current signals, a reduction in the amplitude of approximated versions due to the filtering is not so critical, as it happens to both signals. In contrast, amplitude reduction of frequency components other than 60 Hz may be more significant, since the zeroes of the selected filter can lie on specific harmonic components.

Before introducing the fault location algorithm, the modal components ($\alpha$, $\beta$, and 0) of the voltage signals must be calculated. Its high frequency components are used to identify the fault occurrence. A fault condition is detected when the output of first-level wavelet filter exceeds a stipulated threshold. Two thresholds might be imposed for fault detection. The fault condition may be detected via $\alpha$ and $\beta$-components, and the 0-component indicates whether the ground is involved in the fault or not.

*B. Fault Simulation*

In the following, fault conditions were simulated on a 500 kV three-phase transmission line using ATP (Alternative Transient Program). Fig. 3 shows the simplified diagram of the simulated power system. $Z_{LT}$ denotes the transmission line impedance; $Z_{e1}$ and $Z_{e2}$ are, respectively, the Thévenin equivalent impedances at terminals A and B; $E_{e1}$ and $E_{e2}$ are the Thévenin equivalent voltage source at these terminals.

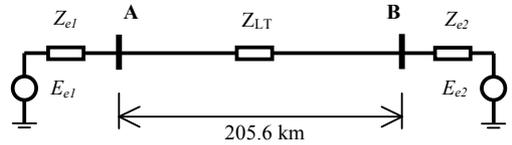

Fig. 3. Simplified Diagram of the Simulated Power System. A three-phase 500 kV transmission line of length 205.6 km was simulated using ATP.

It was considered a totally transposed transmission line with a distributed parameters model (table IV).

TABLE IV
DISTRIBUTED PARAMETERS OF THE TRANSMISSION LINE

| Sequence | Transmission Line Parameters | | |
|---|---|---|---|
| | R (Ω/km) | L (mH/km) | C (nF/km) |
| Positive | 0.0246 | 0.8539 | 13.66 |
| Zero | 0.3818 | 3.732 | 8.61 |

Faults were simulated at 25%, 50%, and 75% of the length of the transmission line, including three single-phase-to-earth faults (A-g, B-g, and C-g); three phase-to-phase faults (AB, AC, and BC); three phase-to-phase-to-earth faults (AB-g, AC-g, and BC-g), and one three-phase fault (ABC). The duration of the ATP-simulation was 8 cycles.

Three fault conditions were simulated: (1) fault occurrences at 4 cycles after beginning simulation, when the voltage of the phase A reaches its maximum; (2) at 4⅛ cycles, and (3) at 4¼ cycles, exactly a zero crossing of the phase A voltage. The action of the protection system opening the line was neglected. Moreover, a null fault resistance was always assumed.

In ATP, a sampling interval of 65.10 μs was selected, but only even indices of the three single-phase voltages and currents at terminal A were saved, yielding 128 samples/cycle.

*C. Fault Analysis with Daubechies and Gegenbauer Filters*

Fault signals were analyzed using 4-coefficient filters: Daubechies (Daub4), Chebyshev (Geg3a1), and Gegenbauer with $\alpha=12$ (Geg3a12). Fig. 4 presents the time-domain scaling and wavelet functions, Fig. 5 the frequency response and Fig. 6 the group delay of the scaling filters.

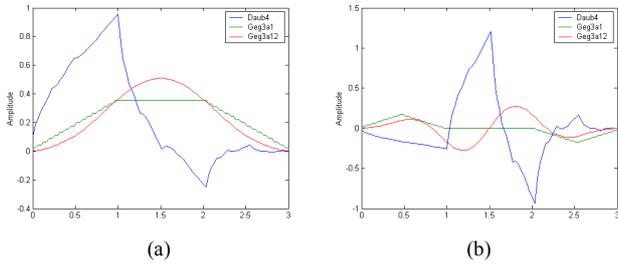

Fig. 4. Analyzing functions derived from 4-coefficient filters after 5 iterations: (a) scaling; (b) wavelet.

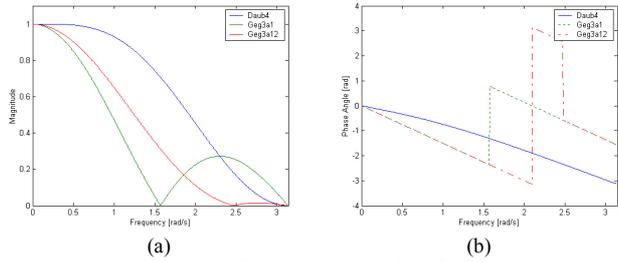

Fig. 5. Frequency response of the analyzing scaling filters: (a) magnitude; (b) phase angle.

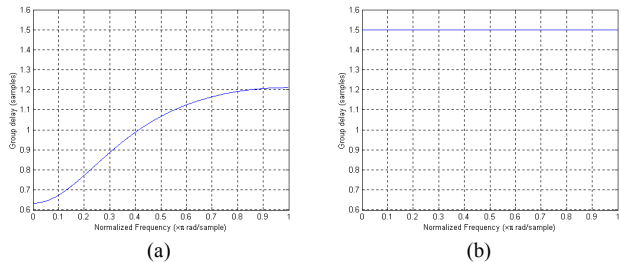

Fig. 6. Group delay of the analyzing scaling filters: (a) Daub4; (b) Geg3a$\alpha$.

The example presented in the sequel illustrates how the analysis can be carried out. Fig. 7 presents the voltage and current waveforms at terminal A for a three-phase fault occurred at 51.4 km from terminal A, at 4⅛ cycles.

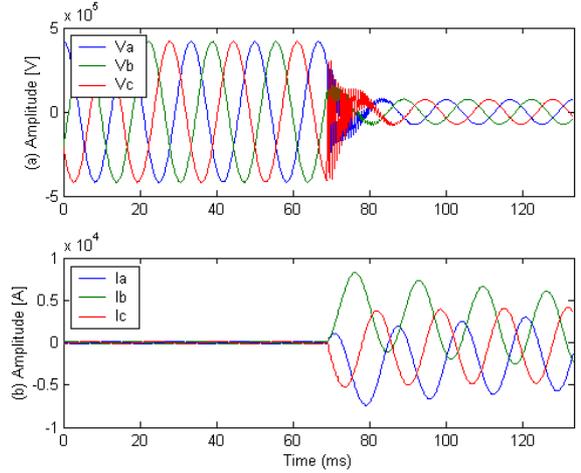

Fig. 7. Signals from a three-phase fault at 51.4 km from Terminal A: (a) voltage signals; (b) current signals.

Details of voltage modal components ($V_\alpha$, $V_\beta$, and $V_0$) of such voltage signals are shown in Fig. 8 by considering only one-stage filter bank decomposition. As expected, only the 0-component was insensitive to the fault indicating that the ground was not involved in the fault. Considering α and β-components, both Daub4 and Geg3a12 filters correctly indicated the fault beginning. In contrast, Geg3a1 identified a misleading fault condition, since it indicated a fault occurrence during the steady-state condition. This can be attributed to its magnitude frequency response (Fig. 5), where the fundamental component appears with a very, however significant, low magnitude. Accordingly, higher *α*-values might be more appropriate to this analysis. The choice of the fault condition threshold depends upon the filter parameters.

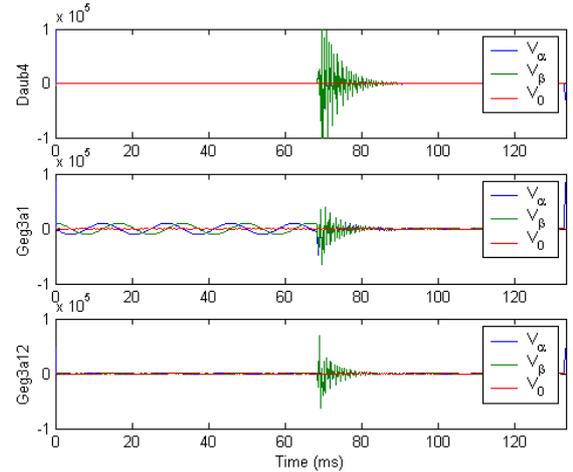

Fig. 8. 1$^{st}$-level detail signal for voltage modal components associated with wavelet filters of: (a) Daub4; (b) Geg3a1; (c) Geg3a12.

Approximated versions of the originals signals are presented in Fig. 9 and Fig. 10 after three-stage filtering decomposition. Comparing these figures with Fig. 6, the removal of high-frequency components can clearly be observed on both voltage and current signals. Furthermore, approximated versions derived from Gegenbauer scaling filters present softer oscillation after the fault incidence than the one shown by Daub4.

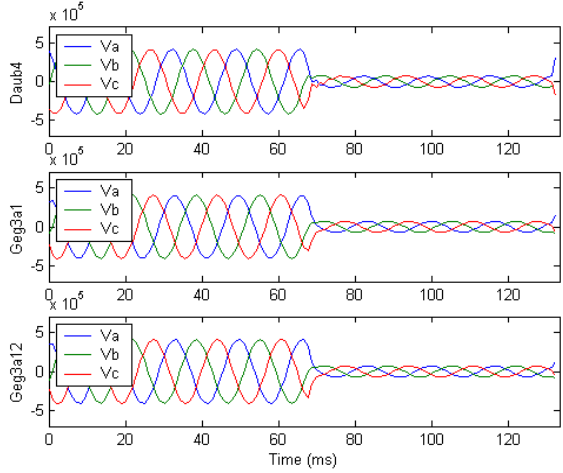

Fig. 9. 3rd-level voltage approximating signal associated with scaling filters of: (a) Daub4; (b) Geg3a1; (c) Geg3a12.

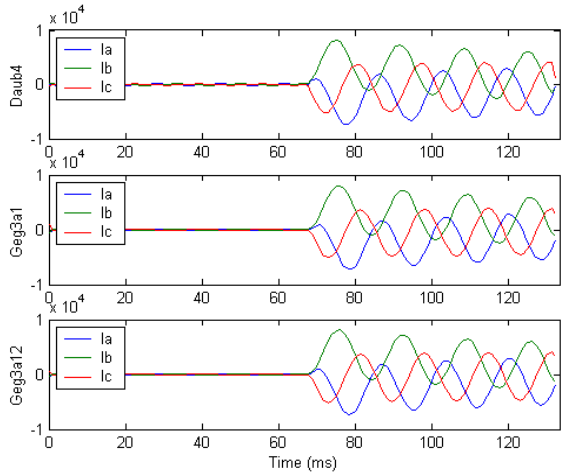

Fig. 10. 3rd-level current approximating signal associated with scaling filters of: (a) Daub4; (b) Geg3a1; (c) Geg3a12.

An estimation of the fault location at each 1-cycle sliding window using these signals is presented in Fig. 11. The detecting algorithm collects one-cycle data window per signal, considering that the oldest sample is disregarded each new entered sample.

The error at the fault location was obtained from:

$$Error = \frac{D_F - D_{FL}}{D_{LT}}, \quad (5)$$

where $D_{FL}$ is the simulated fault distance from terminal A, $D_F$ is the fault distance estimate, and $D_{LT}$ is the length of the transmission line.

The estimate of the fault location is window-dependent probably due to the exponential decaying of the current signals after the fault incidence. Although further requirements to get a more precise indication of the fault location could be incorporated, they were not addressed in this paper, since our primary purpose was to investigate the potential application of Gegenbauer MRA in fault analysis.

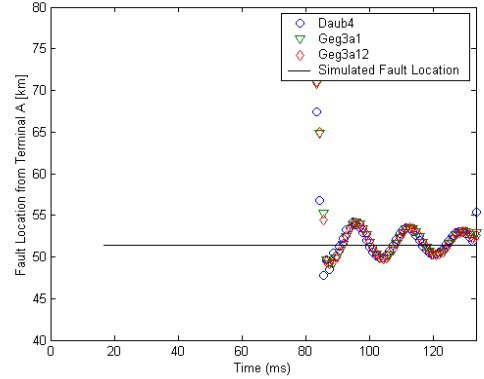

Fig. 11. Fault location estimation from a 1-cycle sliding windowed Fourier transform of the 3rd-level approximating versions associated with scaling filters of Daub4, Geg3a1 and Geg3a12.

The errors of all faults were computed, and for Geg3a1, (Geg3a12), and Daub4 filters, the maximum values were, respectively: 8.90% (8.81%) and 8.05%, for single-phase-to-earth faults; 0.71% (0.78%) and 1.11%, for three-phase faults; 5.47% (5.16%) and 9.32%, for two-phase faults and two-phase-to-earth faults. All values were obtained by reading the 6th window of the eight-cycle monitored data. These figures indicate that both Geg3a1 and Geg3a12 filters can be used in a fault location algorithm.

Table V presents the identifier code (Id) of single-phase-to-earth fault archives derived from the ATP simulations. Ag25_1 describes an A-ground fault at 25% of the length of the transmission line, at position 1, etc. These Id are used in Fig. 12, which shows the fault location error at the 6th-window of the 3rd-level approximating versions associated with scaling filters of Daub4, Geg3a1, and Geg3a12 to all simulated single-phase-to-earth faults.

TABLE V
LABELING SINGLE-PHASE-TO-EARTH FAULT: ID NUMBER VERSUS FAULT TYPE

| Id | Fault Type | Id | Fault Type | Id | Fault Type |
| --- | --- | --- | --- | --- | --- |
| 1 | Ag25_1 | 10 | Ag50_1 | 19 | Ag75_1 |
| 2 | Ag25_2 | 11 | Ag50_2 | 20 | Ag75_2 |
| 3 | Ag25_3 | 12 | Ag50_3 | 21 | Ag75_3 |
| 4 | Bg25_1 | 13 | Bg50_1 | 22 | Bg75_1 |
| 5 | Bg25_2 | 14 | Bg50_2 | 23 | Bg75_2 |
| 6 | Bg25_3 | 15 | Bg50_3 | 24 | Bg75_3 |
| 7 | Cg25_1 | 16 | Cg50_1 | 25 | Cg75_1 |
| 8 | Cg25_2 | 17 | Cg50_2 | 26 | Cg75_2 |
| 9 | Cg25_3 | 18 | Cg50_3 | 27 | Cg75_3 |

The computational effort to implement the algorithm based on Daubechies or Gegenbauer scaling and wavelet filters is commented in the sequel. Due to its symmetry, Gegenbauer filters present a half of the computational effort compared with Daubechies filters, considering same length filters. In particular, setting $\alpha=1$, the $v^{th}$-order Gegenbauer (Chebyshev) filter has identical coefficients, requiring lesser computational effort than any asymmetric filter of same length.

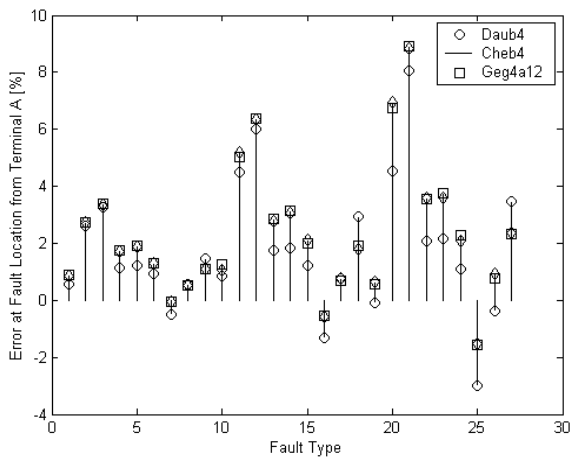

Fig. 12. Fault location error at 6[th]-window of the 3[rd]-level approximating versions associated with scaling filters of Daub4, Geg3a1 and Geg3a12 to single-phase-to-earth faults.

## IV. CONCLUSIONS

The Gegenbauer multiresolution was revisited with the primary purpose of investigating its potential application in fault analysis. Their scaling and wavelet functions have compact support and symmetry, which is attractive from the practical viewpoint. Moreover, Gegenbauer filters introduce the same delay at different frequencies, which may be suitable in specific applications. These filters are offered as an alternative for fault analysis instead of Daubechies filters. The choice of the fault condition thresholds is filter-dependent, and higher $\alpha$-values seem to be more appropriate to fault detection.

This MRA approach can be applied to any subject where the signal reconstruction is not required, supplying more degrees of freedom in the analysis. The choice of $\alpha$ and $v$ of the Gegenbauer filter strongly depend upon what it is intended for.

## VI. BIOGRAPHIES

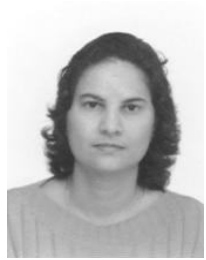

**Luciana R. Soares** was born in Rio de Janeiro, Brazil, on May 30, 1973. She received the B.Sc and M.Sc degrees in electrical engineering from Universidade Federal de Pernambuco (UFPE), Recife, Brazil, in 1997 and 2001. Currently, she is a D.Sc. candidate at the Department of Electronics and Systems at same university (DES-UFPE).

Her special fields of interest include power quality monitoring and analysis, fault analysis, wavelets, and applications of wavelet transforms on power system. She is a student member of CIGRÉ and SBA (Brazilian Society of Automatic).

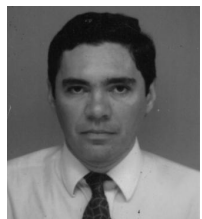

**Hélio M. de Oliveira** (M'1994) was born in Arcoverde, Pernambuco, Brazil, on May 1, 1959. He received the B.Sc. and M.S.E.E. degrees in electrical engineering from Universidade Federal de Pernambuco (UFPE), Recife, Brazil, in 1980 and 1983. Since 1983 he has been with the Department of Electronics and Systems (DES-UFPE). In 1992 he earned the *Docteur* degree from *École Nationale Supérieure des Télécommunications*, Paris.

Dr. de Oliveira was appointed as honored professor by more than 15 undergraduate classes and he authored the book *Análise de Sinais para Engenheiros: Wavelets,* Manole pub., São Paulo, 2004, ISBN 85-204-1624-1, (Signal Analysis for Engineering: Wavelets). His research topics include communication theory, information theory and signal analysis. He is member of IEEE and SBrT (Brazilian Telecommunication Society).